# О ПЕРЕСЕЧЕНИИ ГЛАВНЫХ ИДЕАЛОВ

# В ФУНКЦИОНАЛЬНЫХ ПРОСТРАНСТВАХ

**Меньшикова Э.Б.** E-mail: KhabibullinaE@gmail.com

**Хабибуллин Б.Н.** E-mail: Khabib-Bulat@mail.ru

**Башкирский государственный университет, г. Уфа, РФ**

***Аннотация.*** Даются примеры некоторых колец функций, в которых пересечение двух главных идеалов в этих кольцах не всегда является конечно порождённым идеалом.

***Abstract.*** We give examples of some rings of functions such that the intersection of two principal ideals in these rings is not necessarily a finitely generated ideal.

В работах первого из соавторов [1]-[4] показано, что определённые типы замкнутых идеалов и подмодулей над кольцом многочленов в широких классах топологических колец или модулей над кольцом многочленов, состоящих из голоморфных функций, являются 2-порождёнными. В частности, таковыми оказываются пересечения пар замкнутых главных идеалов/подмодулей в различных весовых кольцах/модулях голоморфных функций. В более общей трактовке 2-порожденными являются в таких кольцах/модулях любые идеалы/подмодули, полностью определяемые их нулевыми множествами, т.е. допускающие локальное описания (определения и подробности см. в [1]-[4]). Для полноты и контраста мы приводим в настоящей заметки простейшие примеры пар главных идеалов в классических пространствах непрерывных действительных функций C[a,b] на отрезке [a,b]$\subset$ (-∞,+∞) и бесконечно дифференцируемых действительных функций $C^{\infty}$(-∞,+∞) , пересечение которых даже не конечно порождённое. Для множества X тождественную функцию на X обозначаем через id: X → X, т.е. id (x) = x для всех x ∈ X. Далее Ideal(*f*) обозначает главный идеал в кольце, порождённый элементом *f* из этого кольца.

**Пример 1.** *Рассмотрим функцию* id *и её модуль* |id| *на* C[-1,1]. *Пусть* Ideal(id) *и* Ideal(|id|) – *главные идеалы в* C[-1,1], *порождённые функциями* id *и* |id|. *Идеал-пересечение* I:=Ideal(id)∩ Ideal(|id|) *не является конечно порождённым.*

В неявной форме этот факт содержится в монографии И. Р. Шафаревича [5; § 4, Пример 8], но мы приведём явное прямое *Доказательство.* Пусть $f$ ∈ I:=Ideal(id)∩Ideal(|id|). Тогда существуют функции

$g, g_1 \in$ C$[-1,1]$, для которых $g \cdot$id = $g_1 \cdot$|id|= $f$. Отсюда $g_1(x) = g(x)$ при всех $0 \le x \le 1$ и $g_1(x) = -g(x)$ при всех $-1 \le x \le 0$. В частности, из непрерывности функции $g_1$ в нуле следует $g$ (0) = – $g$ (0), т.е. $g$ (0) = 0. Таким образом, каждая функция $f$ представляется в виде $f = g \cdot$ id, где $g$ ∈ C[-1,1] и $g$ (0) = 0. Обратно, пусть $f$ — функция указанного вида и положим $g_1(x) = g(x)$ при всех $0 \le x \le 1$ и $g_1(x) = -g(x)$ при всех $-1 \le x \le 0$. Тогда $g_1 \cdot |\mathrm{id}| = g \cdot$id ∈ I . Следовательно, идеал I состоит в точности из функций $f$ вида $g \cdot$ id, где $g$ ∈ C[-1,1] и $g$ (0) = 0, т.е.

I= Ideal(id)∩ Ideal(|id|)={ g · id : $g$ ∈ C[-1,1] и $g$ (0) = 0}.     (1)

Далее воспользуемся методом reductio ad absurdum. Предположим, что существуют натуральное число k и конечный набор функций

$$g_1, g_2, \ldots, g_k \in C[-1, 1], \qquad g_1(0) = g_2(0) = \ldots = g_k\ (0) = 0,$$

для которого набор функций $\{g_1 \cdot \mathrm{id},\ g_2 \cdot \mathrm{id},\ \ldots,\ g_k \cdot \mathrm{id}\}$ порождает идеал I из (1). Рассмотрим функцию

$\mathrm{id} \cdot \left(\max\{\sqrt{|g_1|}, \sqrt{|g_2|}, \ldots, \sqrt{|g_k|}\} + |\mathrm{id}|\right) \in\ \mathrm{I}$ (см. (1)), где $\sqrt{\ } \geq 0$ – арифметический корень. По предположению существует набор непрерывных функций $a_1, a_2, \ldots, a_k \in C[-1, 1]$, с которым имеет место представление

$$\mathrm{id} \cdot \left(\max\{\sqrt{|g_1|}, \sqrt{|g_2|}, \ldots, \sqrt{|g_k|}\} + |\mathrm{id}|\right)$$
$$= a_1 g_1 \cdot \mathrm{id} + a_1 g_1 \cdot \mathrm{id} + \cdots + a_k g_k \cdot \mathrm{id}.$$

Это даёт равенство (в нуле по непрерывности)

$$d \coloneqq \max\{\sqrt{|g_1|}, \sqrt{|g_2|}, \ldots, \sqrt{|g_k|}\} + |\mathrm{id}| = a_1 g_1 + a_1 g_1 + \cdots + a_k g_k.$$

Тогда во всех точках $x \in [-1,0) \cup (0,1]$ справедливы неравенства

$$1 \leq a_1(x) \frac{g_1(x)}{\sqrt{|g_1(x)|} + |x|} + a_2(x) \frac{g_2(x)}{\sqrt{|g_2(x)|} + |x|} + \cdots + a_1(x) \frac{g_k(x)}{\sqrt{|g_k(x)|} + |x|}$$
$$\leq a_1(x)\sqrt{|g_1(x)|} + a_2(x)\sqrt{|g_2(x)|} + \cdots + a_k(x)\sqrt{|g_k(x)|}.$$

При $0 \neq x \to 0$ правая часть здесь стремится к нулю. Противоречие!

**Пример 2.** *Рассмотрим функцию в* $C^{\infty}(-\infty,+\infty)$ *соответственно нечётную и чётную бесконечно дифференцируемые функции*

$$c(x) = \begin{cases} e^{-\frac{1}{x}} \text{ при } x \in (0, +\infty), \\ \quad 0 \quad \textit{при}\ x = 0, \\ -e^{\frac{1}{x}} \ \textit{при}\ x \in (-\infty, 0), \end{cases} \qquad c_0 \coloneqq |c|.$$

*Идеал-пересечение* I:= Ideal($c$)∩ Ideal($c_0$) *не является конечно порождённым.*

Доказательство опускаем, поскольку оно аналогично доказательству предыдущего Примера 1 и во многом сводится ко второй части [5; § 4, Пример 8].

**Замечание.** Как отмечено в [5; § 4, после Примера 8], любые замкнутые идеалы в топологических кольцах непрерывных C[a,b] и бесконечно дифференцируемых $C^{\infty}(-\infty,+\infty)$ функций с их естественными топологиями являются главными и аналоги примеров типа 1 и 2 для замкнутых идеалов невозможны.

Отметим здесь, что построение аналогов Примеров 1 и 2 в различных кольцах голоморфных функций в областях из комплексной плоскости представляет интерес и, по-видимому, такие конструкции будут значительно более сложными. Некоторым подспорьем для построения таких примеров могут быть статьи Л.И. Гуревича [6] и первого из соавторов [7], где рассмотрены контрпримеры к задаче спектрального анализа/синтеза, зачастую двойственной к проблеме локального описания идеалов/подмодулей в кольцах/модулях голоморфных функций. Эти вопросы в отношении голоморфных функций предполагается рассмотреть в ином месте.

## ЛИТЕРАТУРА